\documentclass[a4paper]{eccomas_2022}
\usepackage{graphicx}
\usepackage{amsmath}
\usepackage{amsfonts}
\usepackage{amssymb}

\usepackage{graphicx}
\usepackage{amsmath}
\usepackage{amsfonts}
\usepackage{amssymb}
\usepackage{algorithm,todonotes}
\usepackage[noend]{algpseudocode}
\usepackage{booktabs}
\DeclareMathOperator*{\argmax}{arg\,max}
\DeclareMathOperator*{\argmin}{arg\,min}
\newcommand{\ep}{\varepsilon_{\mathcal{T}}}

\usepackage{mathtools}
\DeclarePairedDelimiter\floor{\lfloor}{\rfloor}
\newcommand{\R}{\mathbb R}

\title{Feedback reconstruction techniques for optimal control problems on a tree structure}

\author{ALESSANDRO ALLA$^1$ AND  LUCA SALUZZI$^2$}

\heading{A. Alla, M. Falcone and L. Saluzzi}

\address{$^{1}$ Dipartimento di Scienze molecolari e nanosistemi\\
Università Ca' Foscari Venezia\\\
e-mail: alessandro.alla@unive.it
\and
$^{2}$ Department of Mathematics\\ Imperial College London\\
email: l.saluzzi@imperial.ac.uk
}


\abstract{The computation of feedback control using Dynamic Programming equation is a difficult task due the {\em curse of dimensionality}. The tree structure algorithm is one the methods introduced recently that mitigate this problem. The method computes
the value function avoiding the construction of a space grid and the need for interpolation
techniques using a discrete set of controls. However,  the computation of the control is strictly linked to control set chosen in the computation of the tree. Here, we extend and complete the method selecting a finer control set in the computation of the feedback. This requires to use an interpolation method for scattered data which allows us to reconstruct the value function for nodes not belonging to the tree.
The effectiveness of the method is shown via a numerical example.}

\begin{document}

\section{Introduction}

The computation of feedback control for differential equations is an important topic due to applications in real life problems. Usually, one uses the dynamic programming principle and the Hamilton--Jacobi--Bellman (HJB) equations to derive the control in feedback form (see e.g. \cite{BCD97} for a complete description of the method). The major issue of this approach is that the solution of the HJB equation is not analytical and we need to build numerical approximations. Although there exists a huge literature on the approximation (see e.g. \cite{FF13}),  numerical methods suffer from the curse of dimensionality, namely the complexity of the problem increases as the dimension of the system we want to control does. In the last decades, there were a tremendous effort in mitigating the curse of dimensionality using different methods such as: model order reduction \cite{KVX04, AFV17}, spectral methods \cite{KK18}, max-plus algebra \cite{ME07, ME09}, neural networks \cite{DLM20, DM21}, tensor decomposition \cite{DKK21, OSS22, DKL22}, sparse grids method \cite{GK17} and radial basis functions \cite{AOS21}.

Recently, in \cite{AFS18} it has been introduced a tree structure algorithm to approximate the HJB equation for finite horizon problem. The value function is computed using a DP algorithm on a tree structure algorithm (TSA) constructed by the time discrete dynamics. In this way there is no need to build a fixed space triangulation and to project on it: the tree will guarantee a perfect matching with the discrete dynamics and drop off the cost of the space interpolation allowing for the solution of very high-dimensional problems. Moreover, a pruning technique has been implemented to reduce the number of branches and the exponential complexity of the tree. Error estimates have been derived in \cite{SAF19} for the TSA, including the pruning technique, to guarantee first order convergence. Later, the method has been extended to high order methods in \cite{AFS18b}, state constraint problems \cite{AFS20} and coupled with model order reduction to deal with large scale problems \cite{AS19}.

In this work, we conclude the study of the TSA explaining how to build the feedback control. Indeed, in the works presented before the control was linked to the nodes of the tree and it was not able to be obtained for different initial conditions. In fact, the methods could not reconstruct the control for points outside of the tree nodes. Here, we propose two algorithms based on scattered interpolation to overcome this limit. In particular, given the tree structure and the information of the value function on it, one may apply interpolation operators on scattered data to construct the value function, and hence the feedback map, on points not belonging to the tree.

The outline of the paper is the following. In Section 2, we recall the tree structure algorithm and the main ingredients for the resolution of the dynamic programming principle on the TSA. Section 3  is devoted to the study of our new method to reconstruct feedback control based on the use of interpolation techniques on scattered dataset. Finally, in Section 4 we present a numerical test to show the effectiveness of the proposed methodology.

\section{The tree structure algorithm}

In this section we will recall the {\em finite horizon control problem} and its approximation by the TSA (see \cite{AFS18} for a complete description of the method). Let us consider the following dynamics
\begin{equation}\label{eq}
\left\{ \begin{array}{l}
\dot{y}(s)=f(y(s),u(s),s), \;\; s\in(t,T],\\
y(t)=x\in\R^d.
\end{array} \right.
\end{equation}
where $y:[t,T]\rightarrow\R^d$ is the solution, $u:[t,T]\rightarrow\R^m$ is the control, $f:\R^d\times\R^m\times[t,T]\rightarrow\R^d$ is the dynamics and
\[\mathcal{U}=\{u:[t,T]\rightarrow U, \mbox{measurable} \}
\]
is the set of admissible controls within the compact set $U\subset \R^m$. 
We assume that there exists a unique solution for \eqref{eq} for each $u\in\mathcal{U}$. 
The cost functional we want to minimize reads

\begin{equation}\label{cost}
 J_{x,t}(y,u):=\int_t^T L(y(s),u(s),s)e^{-\lambda (s-t)}\, ds+g(y(T))e^{-\lambda (T-t)},
\end{equation}
where $L:\R^d\times\R^m\times [t,T]\rightarrow\R$ is the running cost and $\lambda\geq0$ is the discount factor. 
Finally, the optimal control problem is
\begin{equation}\label{ocp}
\min_{u\in \mathcal{U}} J_{x,t}(y,u), \mbox{ subject to } y(\cdot;u) \mbox{ solution of  \eqref{eq}}
\end{equation}
We assume that the functions $f, L$ and $g$ are bounded and Lipschitz-continuous with respect to the first variable to guarantee existence and uniqueness of the control problem \eqref{ocp}.

The value function is defined as follows
\begin{equation}
v(x,t):=\inf\limits_{u\in\mathcal{U}} J_{x,t}(u)
\label{value_fun}
\end{equation}
and satisfies the DPP, i.e. for every $\tau\in [t,T]$:
\begin{equation}\label{dpp}
v(x,t)=\inf_{u\in\mathcal{U}}\left\{\int_t^\tau L(y(s),u(s),s) e^{-\lambda (s-t)}ds+ v(y(\tau),\tau) e^{-\lambda (\tau-t)}\right\}.
\end{equation}
From \eqref{dpp}, one can derive the HJB equation for every $x\in\R^d$, $s\in [t,T)$: 
\begin{equation}\label{HJB}
\left\{
\begin{array}{ll} 
-\dfrac{\partial v}{\partial s}(x,s) +\lambda v(x,s)+\max\limits_{u\in U }\left\{-L(x, u,s)- \nabla v(x,s) \cdot f(x,u,s)\right\} = 0, \\
v(x,T) = g(x).
\end{array}
\right.
\end{equation}
Finally, the computation of the feedback control is straightforward, assuming the value function is known:
\begin{equation}\label{feedback}
u^*(x):=  \argmax_{u\in U }\left\{-L(x,u,t)- \nabla v(x,t) \cdot f(x,u,t)\right\}. 
\end{equation}

Since equation \eqref{HJB} is a first non-linear PDE, it is hard to find an exact solution and numerical algorithms should take into account discontinuities in the gradient (see \cite{FF13} and the references therein). Introduced a time discretization of \eqref{HJB} with a time step $\Delta t: = [(T-t)/\overline N]$ and $\overline{N}$ number of steps, it is possible to consider the discrete version of the DPP \eqref{dpp}. More precisely, for $n= \overline{N}-1,\dots, 0$ and every $x\in \R^d$ we have 

\begin{equation}\label{SL}
V^{n}(x)=\min\limits_{u\in U}[\Delta t\,L(x, u, t_n)+e^{-\lambda \Delta t}V^{n+1}(x+\Delta t f(x, u, t_n))], 
\end{equation}
where $t_n=t+n \Delta t,\, t_{\overline N} = T$, and $V^n(x):=V(x, t_n).$ The iterative scheme \eqref{SL} is coupled with the terminal condition 
\begin{equation}\label{SL_T}
V^{\overline{N}}(x)=g(x).
\end{equation}
In \eqref{SL} we use an explicit Euler scheme for a first order approximation to simplify the presentation (the high-order extension has been presented in \cite{AFS18b}). The term $V^{n+1}(x+\Delta t f(x, u, t_n))$ is usually obtained via interpolation on a fixed grid since $x+\Delta t f(x, u, t_n)$ is not a grid point (see \cite{FF13} for more details on this step). To bypass the interpolation step a tree structure is built where all the possible combinations of the term $x+\Delta t f(x, u, t_n)$ are computed  for different values of $u$. \\
First of all, let us consider a discrete version of the control domain, say $U=\{u_1,...,u_M \}$ with $M$ controls.
 We will denote the tree by $\mathcal{T}:=\cup_{j=0}^{\overline{N}} \mathcal{T}^j,$ where each $\mathcal{T}^j$ contains the nodes of the tree at time $t_j$. The first level $\mathcal{T}^0 = \{x\}$ is simply formed by the initial condition $x$. Starting from the initial condition $x$, we discretize the dynamics using  e.g. an explicit Euler scheme and we consider all the nodes obtained with different discrete controls $u_i \in U $
$$\zeta_{i}^1 = x+ \Delta t \, f(x,u_i,t_0),\qquad i=1,\ldots,M.$$ 
Therefore, we have $\mathcal{T}^1 =\{\zeta_1^1,\ldots, \zeta^1_M\}$.
The procedure can be easily iterated for each node of the level, obtaining at time $t_n$ the level $\mathcal{T}^n$:
$$\mathcal{T}^n = \{ \zeta^{n-1}_i + \Delta t f(\zeta^{n-1}_i, u_j,t_{n-1}) \}_{j=1}^{M}\quad i = 1,\ldots, M^{n-1}.$$ 
The entire tree can be represented in short as
 $$\mathcal{T}:= \{ \zeta_{i}^n  \}_{i=1}^{M^n},\quad n=0,\ldots \overline{N},$$ 
where the nodes $\zeta^n_i$ are the results of the dynamics at time $t_n$ with the controls $\{u_{j_k}\}_{k=0}^{n-1}$:
\begin{equation*}
\zeta_{i_n}^n = \zeta_{i_{n-1}}^{n-1} + \Delta t f(\zeta_{i_{n-1}}^{n-1}, u_{j_{n-1}},t_{n-1})
= x+ \Delta t \sum_{k=0}^{n-1} f(\zeta^k_{i_k}, u_{j_k},t_k), 
\end{equation*}
with $\zeta^0 = x$, $i_k = \floor*{\dfrac{i_{k+1}}{M}}$ and $j_k\equiv i_{k+1} \mbox{mod } M$.

Although it is possible to deal with arbitrary high-dimensional problems, the construction of tree may be expensive since  $|\mathcal{T}|= O(M^{\overline{N}})$,
where $M$ is the number of discrete controls and $\overline{N}$ is the number of time steps. This leads to an exponential growth of the cardinality and it may be infeasible to apply the algorithm due to the huge amount of memory allocations, when $M$ or $\overline{N}$ are too large.
To mitigate this exponential growth a pruning criteria has been introduced.
Defining a threshold $\ep>0$, several branches of the tree can be cut off according to the distance between nodes. More precisely, if two nodes satisfies the following criteria
\begin{equation}\label{tol_cri}
\begin{array}{cc}
\Vert \zeta^n_i-\zeta^n_j \Vert \le \ep,
\mbox{ for  }i\ne j \mbox{ and } n = 0,\ldots, \overline{N}.
\end{array}
\end{equation}
they can be merge, leading to a great gain in memory storage and computational time, keeping the same order of convergence (\cite{SAF19}). 


The computation of the numerical value function $V(x,t)$ will be performed on the tree nodes in space as 
\begin{equation}\label{num:vf}
V(x,t_n)=V^n(x), \quad \forall x \in \mathcal{T}^n, 
\end{equation}

where $t_n=t+ n \Delta t$, following directly from the DPP. The tree $\mathcal{T}=\cup_{j=0}^{\overline{N}} \mathcal{T}^j$ given by the TSA defines a grid and we can write a time discretization on it for \eqref{HJB} as follows: 

\begin{eqnarray}\label{HJBt2}
\begin{cases}
V^{n}(\zeta^n_i)= \min\limits_{u\in U} \{e^{-\lambda \Delta t} V^{n+1}(\zeta^n_i+\Delta t f(\zeta^n_i,u,t_n)) +\Delta t \, L(\zeta^n_i,u,t_n) \}, \hbox{ for }\zeta^n_i \in \mathcal{T}^n, n \le \overline{N}-1,\nonumber\\
V^{\overline{N}}(\zeta^{\overline{N}}_i)= g(\zeta_i^{\overline{N}}),  \qquad\quad \hbox { for }   \zeta_i^{\overline{N}} \in \mathcal{T}^{\overline{N}}.\qquad
\end{cases}
\end{eqnarray}

Since the set of controls $U$ is discrete, the minimization will be computed by comparison.
\medskip
\subsection{TSA and Model Order Reduction}
\label{POD}

The Tree Structure Algorithm has been also coupled with Proper Orthogonal Decomposition (see e.g. \cite{V13})  techniques in \cite{AS19}. The idea is to consider the projection of the dynamics onto a subspace spanned by particular orthogonal basis function in order to reduce the dimension and the complexity of the problem. Here, we sketch briefly the main concepts, the interested reader will find more details in \cite{AS19}.

Given the full dimensional dynamics \eqref{eq}, the TSA can be applied with few time steps and few controls to explore the manifold of all possible solutions. Collected all the nodes of the tree in a matrix $Y \in \mathbb{R}^{d \times N}$, called the \emph{snapshots matrix}, one can operate a Singular Value Decomposition of the matrix $Y$, obtaining $Y=\Psi \Sigma V^\top$, where $\Psi \in \mathbb{R}^{d \times d}$ and $V \in \mathbb{R}^{N \times N}$ are orthogonal matrices and $\Sigma \in \mathbb{R}^{d \times N}$ is a diagonal matrix with diagonal entries $\{\sigma_i\}_{i=1}^{\min\{d,N\}}$. The first $\ell \ll \min\{d,N\}$ columns of the matrix $\Psi$, $i.e.$ $\Psi^\ell =\{\psi_1, \ldots, \psi_\ell\}$, will represent the solution of the following minimization problem
\begin{equation*}\label{pbmin}
\min_{ {{\psi}}_1,\ldots,{{\psi}}_\ell\in\R^d} \sum_{j=1}^N \left|{ y}^j-\sum_{i=1}^\ell \langle { y}^j,{{\psi}}_i\rangle{{\psi}}_i\right|^2\quad \mbox{such that }\langle {{\psi}}_i,{{\psi}}_j\rangle=\delta_{ij}.
\end{equation*}
The number of basis $\ell$ can be fixed according to a criterium related to the projection error. More precisely, given a tolerance $\tau \in [0,1]$ the parameter $\ell$ can be chosen such that
\begin{equation}\label{POD_cri}
\mathcal{E}(\ell)=\dfrac{\sum_{i=\ell+1}^{\min\{d,N\}} \sigma_i^2}{\sum_{i=1}^{\min\{d,N\}} \sigma_i^2} \le \tau.
\end{equation}
The constructed matrix $\Psi^\ell$ will be then employed in the projection of the dynamics \eqref{eq}. Indeed, assuming the ansatz $y(t) \approx \Psi^\ell y^\ell(t)$ with $y^\ell(t) \in \mathbb{R}^\ell$, and the orthogonality of $\Psi^\ell$ the reduced dynamics reads
\begin{equation*}\label{pod_sys}
\left\{\begin{array}{l}
\dot{ y}^\ell(s)= (\Psi^\ell)^\top f( \Psi^\ell  y^\ell(s),u(s),s),\\
{ y}^\ell(t)={ \Psi^\ell}^\top x.
\end{array}\right.
\end{equation*}
In our numerical experiment we will use this method for the control of the heat equation to reduce the dimensionality of the problem.

\section{Feedback reconstruction and closed-loop control}
\label{subsec}

In this section we are going to present a technique to retrieve the feedback control based on the knowledge of the value function on the nodes of the tree. We are going to introduce two possible post-processing reconstruction which will be compared in the section of the numerical test.

The computational cost for the construction of the full tree is exponential in the number of discrete controls, for this reason it is better to consider few controls for the tree construction and for the resolution of the HJB equation. Once the value function is obtained on a tree-structure, we consider a post-processing procedure which takes into account a finer control set. This is possible thanks to the formula for the synthesis of the feedback control
\begin{equation} \label{feed:tree2}
u_{n}^{*}:=\argmin\limits_{u\in \tilde{U}} \left\{ e^{-\lambda \Delta t}V^{n+1}(x+\Delta t f(x,u,t_n)) +\Delta t \, L(x,u,t_n) \right\},
\end{equation}
where the $argmin$ will be computed on a finer set $\tilde{U}$ with respect to the initial set $U$. This minimization can be computed again by comparison, but we need to introduce an interpolation step on scattered data. In low dimension (i.e. two or three) one may consider a Delaunay triangulation of the data and then perform an interpolation on the triangulation. In high dimensions, the triangulation becomes infeasible and one has to proceed in different ways, for example via e.g. kernel methods (\cite{F15}) or via Model Order Reduction Techniques (\cite{AS19}). In the numerical test we will consider the latter case, considering a POD reduction of the dynamical system.

The method for the feedback reconstruction based on a minimization by comparison on a finer control set is presented in Algorithm \ref{alg_4}. 
\begin{algorithm}[H]
\caption{Feedback reconstruction via comparison on a finer control set}
\label{alg_4}
\begin{algorithmic}[1]
 \State{Computation of the tree $\mathcal{T}$ and value function $\{V^k\}_k$ with control set $U$.}
 \State{Fix a new control set $\tilde{U}\supset U$ and $\zeta^0_*=x.$}
\For{$n=0,...,\overline{N}-1$}
 \For{$u_j \in \tilde{U}$}
    \State{$\zeta_j=\zeta^{n}_*+ \Delta t f(\zeta^{n}_*,u_j,t_{n})$}
    \State{Compute $V(\zeta_j,t_{n+1})$ via scattered interpolation with $(\mathcal{T}^{n+1},V^{n+1})$}
    \EndFor
    \State{$u_{n}^{*}:=\argmin\limits_{u_j\in \tilde{U}} \left\{ e^{-\lambda \Delta t}V(\zeta_j,t_{n+1}) +\Delta t \, L(\zeta^n_*,u_j,t_n) \right\}$}
    \State{$\zeta^{n+1}_*=\zeta^{n}_*+ \Delta t f(\zeta^{n}_*,u^*_n,t_{n})$}
    \EndFor
\end{algorithmic}
\end{algorithm}
The interpolation on scattered data can be computed via the MATLAB function

\texttt{scatteredInterpolant}. If the dynamics $f$ is autonomous, by Remark $3.3$ in \cite{AFS18} we know that we can compute at time $t_n$ the value function on the sub-tree $\cup_{k=0}^n \mathcal{T}^k$. In this case, in step $6$ we can compute the scattered interpolation with $\cup_{k=0}^n(\mathcal{T}^k,V^k)$, guaranteeing more information for a more efficient interpolation. 

Now, let us consider a dynamics $f$ affine in the control $u\in \mathbb{R}$. By Remark $3.1$ in \cite{AFS18} we know that all the tree sons of a node lay on a segment. In this case we can apply one dimensional interpolation, for instance quadratic interpolation if we consider three discrete controls for each iteration. The quadratic interpolation is a reasonable choice in the Linear Quadratic Regulator case since we know that the value function is quadratic and then we do not introduce interpolation error in this case. Moreover, let us suppose that the running cost $L$ is of the form $L(x,u,t)=g(x,t)+\gamma|u|^2+\delta u$. In Algorithm \ref{alg_5} we describe this procedure based on a quadratic interpolation, fixing $\lambda=0$ for simplicity. In step $9$ the operator $\mathcal{P}_U$ stands for the projection operator onto the set $U$. We will use and compare these two techniques in the numerical experiment, where we will consider the optimal control for the heat equation.

\begin{algorithm}[H]
\caption{Feedback reconstruction via quadratic interpolation}
\label{alg_5}
\begin{algorithmic}[1]
\State{Computation of the tree $\mathcal{T}$ and value function $\{V^k\}_k$ with control set $U=\{u_1,u_2,u_3 \}$.}
 \State{$\zeta^0_*=x.$}
\For{$n=0,...,\overline{N}-1$}
 \For{$u_j \in U$}
    \State{$\zeta(u_j)=\zeta^{n}_*+ \Delta t f(\zeta^{n}_*,u_j,t_{n})$}
    \State{Compute $V(\zeta(u_j),t_{n+1})$ via scattered interpolation with $(\mathcal{T}^{n+1},V^{n+1})$}
    \EndFor
    \State{$V(\zeta(u),t_{n+1}) \approx au^2+bu+c, \quad \forall u\in [u_1,u_3]$}
    \If{$a+\Delta t\, \gamma>0$ }
    \State{$u_{n}^{*}=\mathcal{P}_{[u_1,u_3]}\left(-\frac{b+\Delta t \, \delta}{2(a+\Delta t\, \gamma)}\right) $}
    \Else
    
    \State{$u_{n}^{*}=\argmin\limits_{u_i \in \{u_1,u_3 \} } \left\{ V(\zeta(u_i),t_{n+1}) +\Delta t \, L(\zeta^n_*,u_i,t_n) \right\}$}
    
\EndIf
    \State{$\zeta^{n+1}_*=\zeta^{n}_*+ \Delta t f(\zeta^{n}_*,u^*_n,t_{n})$}
    \EndFor
\end{algorithmic}
\end{algorithm}
\section{Numerical experiment: Heat Equation}

\begin{figure}[htbp]
\centering
       \includegraphics[scale=0.45]{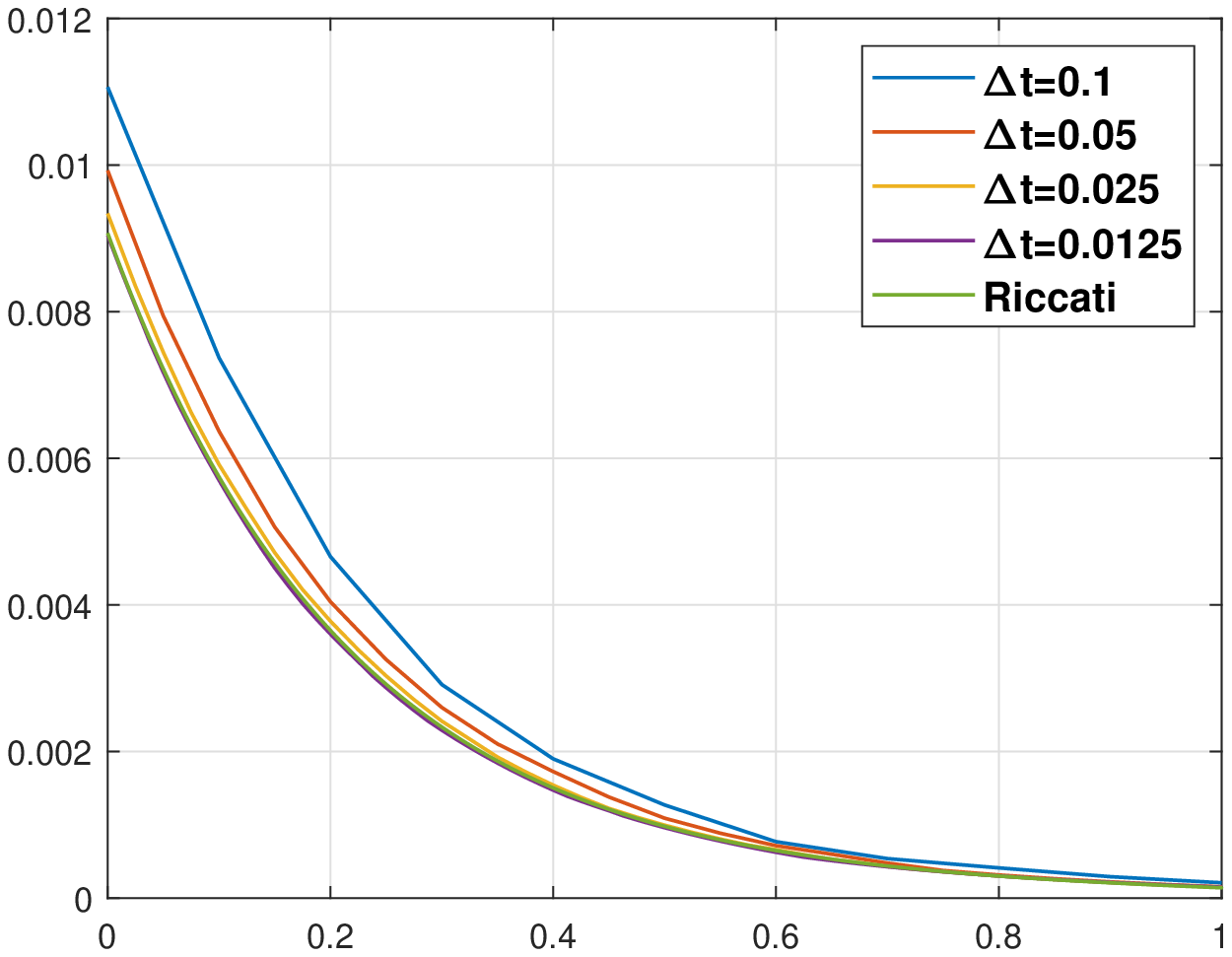}
       \includegraphics[scale=0.45]{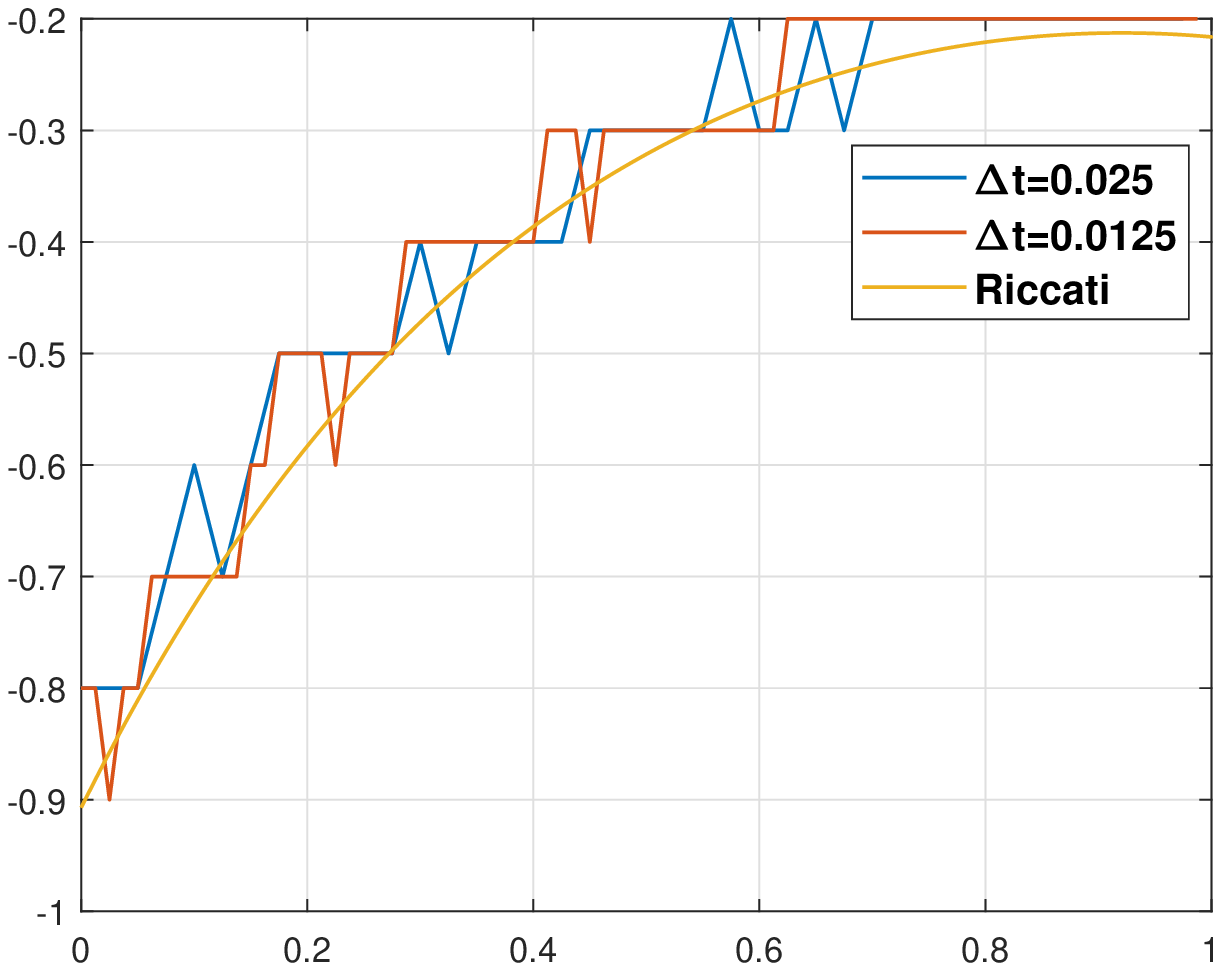} 
       \caption{Test 1: Cost functional (left) and optimal control (right) with $11$ discrete controls. }
       \label{fig3:cost11}
	\end{figure}

In this example we consider the one dimensional heat equation with homogeneous Dirichlet boundary conditions
\begin{equation}
\begin{cases}
\partial_t y(x,t)= \sigma  y_{xx}(x,t) +  y_0(x)u(t) & (x,t) \in \Omega \times [0,T],
\\ y(x,t) =0 & (x,t) \in \partial \Omega \times [0,T], \\
y(x,0)=y_0(x) & x \in \Omega,
\end{cases}
\label{pde1}
\end{equation}

where the control $u(t)$ is taken in the admissible set $\mathcal{U}=\{u:[0,T]\rightarrow [-1,0] \}$, $\sigma=0.15$, $T=1$ and $\Omega=[0,1]$. The space domain is discretized with $1000$ nodes, leading to a dynamical system in dimension $\mathbb{R}^{1000}$. Since the problem is truly high dimensional, we are going to apply the POD reduction introduced in Section \ref{POD}.
We create a rough tree with two discrete controls and $\Delta t=0.1$ and we fix the tolerance $\tau=10^{-4}$ in \eqref{POD_cri}. This procedure yields a reduced dimension $\ell=2$. In this case two POD basis get enough information for a quasi-complete description of the system. Then, given the reduced dynamics, we compute the reduced tree structure and the value function with $11$ discrete controls and we compare the numerical value function computed by TSA-POD and by the Riccati equation (see e.g. \cite{AFIJ03}) according to the following errors
\begin{equation}
Err_{2}= \frac{\sum_{n=0}^N |V(y^n_*,t_n)-v(y^n_{R},t_n)|^2}{\sum_{n=0}^N |v(y^n_{R},t_n)|^2}, \quad
 Err_{\infty}=   \frac{ \max\limits_{n=0,...,N} |V(y^n_*,t_n)-v(y^n_{R},t_n)|}{  \max\limits_{n=0,...,N} |v(y^n_{R},t_n)|},
 \label{errors}
\end{equation}

where $\{y^n_*\}_n$ is the optimal trajectory computed via POD-TSA, while $\{y^n_{R}\}_n$ is the solution of the time dependent Riccati equation. The results are presented in Table \ref{test41:11contr}. The column Pruned/Full refers to the ratio between the cardinalities of the pruned tree and the full tree. Moreover, in Table \ref{test1:11contr} we recall the values obtained for the full dimensional problem. The coupling of the TSA with POD leads to more accurate results in less time. In particular, we see that for $\Delta t=0.0125$ we obtain a speed-up of a factor $7$ and a reduction of order $4$ for the cardinality of tree. In Figure \ref{fig3:cost11} we can observe the convergence of the cost functional and the approximation of the optimal control.

\begin{table}[htbp]
		\centering
		\begin{tabular}{cccccccc}
					\toprule
					 $\Delta t$ & Nodes & Pruned/Full & CPU & $Err_{2}$ &  $Err_{\infty} $ & $Order_{2}$  & $Order_{\infty}$  \\
					\midrule
					0.1 &     134 & 4.3e-10	 &  \phantom{x}0.1s&    0.244 &   0.220 &   & \\
					0.05 &   825 &  1.0e-19	 &  \phantom{x} 0.56s &   0.102 &  9.4e-2 &  1.25 &  1.22 \\
					0.025 &   11524 & 2.1e-39	 &  \phantom{x} 8.74s & 3.1e-2 & 3.0e-2  &   1.73 &    1.67 \\
					0.0125 &  194426 & 7.8e-80 &  \phantom{x}151s & 1.0e-2   &   8.2e-3 &  1.60 &   1.85  \\
			\bottomrule		
		\end{tabular}
   	\caption{Test 1: Error analysis and order of convergence for TSA-POD method with $\ep= \Delta t^2$, $11$ discrete controls and $2$ POD basis.}
	\label{test41:11contr}
	\end{table}
	
\begin{table}[htbp]
		\centering
		\begin{tabular}{cccccccc}
					\toprule
					 $\Delta t$ & Nodes & Pruned/Full & CPU & $Err_{2}$ &  $Err_{\infty} $ & $Order_{2}$  & $Order_{\infty}$  \\
					\midrule
					0.1 &     134 & 4.7e-09	 &  \phantom{x}0.14s&    0.279 &   0.241 &   & \\
					0.05 &   863 & 1.2e-18	 &  \phantom{x} 0.65s &   0.144 &  0.118 &  0.95 &  1.03 \\
					0.025 &   15453 & 3.1e-38	 &  \phantom{x} 12.88s & 5.5e-2  & 5.3e-2  &  1.40 &    1.17 \\
					0.0125 &  849717 & 3.8e-78 &  \phantom{x}1.1e3s & 1.6e-2   &   1.6e-2 &  1.77 &  1.42  \\
			\bottomrule		
		\end{tabular}
   	\caption{Test 1: Error analysis and order of convergence for forward Euler scheme of the TSA with $\ep= \Delta t^2$ and $11$ discrete controls.}
	\label{test1:11contr}
	\end{table}
	
\paragraph{Feedback reconstruction}

In this paragraph we are going to apply the feedback reconstruction techniques introduced in Section \ref{subsec}. We consider again $2$ POD basis and we solve the optimal control problem via POD-TSA with $3$ discrete controls. The results for this case are presented in Table \ref{test1:2contr}.

\begin{table}[htbp]
		\centering
		\begin{tabular}{cccccccc}
					\toprule
					 $\Delta t$ & Nodes & Pruned/Full & CPU & $Err_{2}$ &  $Err_{\infty} $ & $Order_{2}$  & $Order_{\infty}$  \\
					\midrule
					0.1 &     122 & 4.6e-04	 &  \phantom{x}  0.02s&    0.376 &   0.283 &   & \\
					0.05 &   689 &  4.4e-08	 &  \phantom{x} 0.19s &   0.178 &  0.136 &  1.08 &  1.06 \\
					0.025 &   9536 & 1.7e-16	 &  \phantom{x} 2.3s & 0.107 & 6.9e-2  &   0.73
 &    0.98
 \\
					0.0125 &  155293 & 2.3e-34 &  \phantom{x}37s & 0.0655   &   0.0394 &  0.71 &   0.80  \\
			\bottomrule		
		\end{tabular}
   	\caption{Test 1: Error analysis and order of convergence for TSA-POD method with $\ep= \Delta t^2$, $3$ discrete controls and $2$ POD basis.}
	\label{test1:2contr}
	\end{table}
	Then, we pass to the post-processing procedure: we consider Algorithm \ref{alg_4}, computing the optimal trajectory/control with a finer control set. The control set $\tilde{U}$ has now $100$ discrete controls. In Table \ref{test1:2contrcomp} we present the errors and the orders according to the definition \eqref{errors}.
	
	\begin{table}[htbp]
		\centering
		\begin{tabular}{cccccccc}
					\toprule
					 $\Delta t$ &  CPU & $Err_{2}$ &  $Err_{\infty} $ & $Order_{2}$  & $Order_{\infty}$  \\
					\midrule
					0.1       &  \phantom{x} 0.03s&    0.315 &   0.250 &   & \\
					0.05     	 &  \phantom{x}  0.07s &   9.6e-2 &  0.100 & 1.71 &   1.32 \\
					0.025   	 &  \phantom{x}0.74s & 2.5e-2 & 3.1e-2  &  1.93 &    1.68 \\
					0.0125  &  \phantom{x} 25s & 1.4e-2   &   9.0e-3 &  0.89 &   1.81  \\
			\bottomrule		
		\end{tabular}
   	\caption{Test 1: Error analysis and order of convergence for TSA-POD method with $\ep= \Delta t^2$, $2$ POD basis and reconstruction with $100$ controls.}
	\label{test1:2contrcomp}
	\end{table}
	
Since we have constructed the tree with $3$ discrete controls, we can apply also the quadratic feedback reconstruction presented in Algorithm \ref{alg_5}. The results for this case are presented in Table \ref{test1:2contrpar}. It is possible to notice that the CPU time and the errors are similar. Tables \ref{test1:2contrcomp}-\ref{test1:2contrpar}  can be now compared with Table \ref{test41:11contr} in which we were solving the optimal control problem without introducing the feedback reconstruction technique. The numerical errors in both norms are comparable in the three cases studied, while in term of computational cost we can notice a speed-up of almost 3 orders.

\begin{figure}[htbp]
\centering
       \includegraphics[scale=0.45]{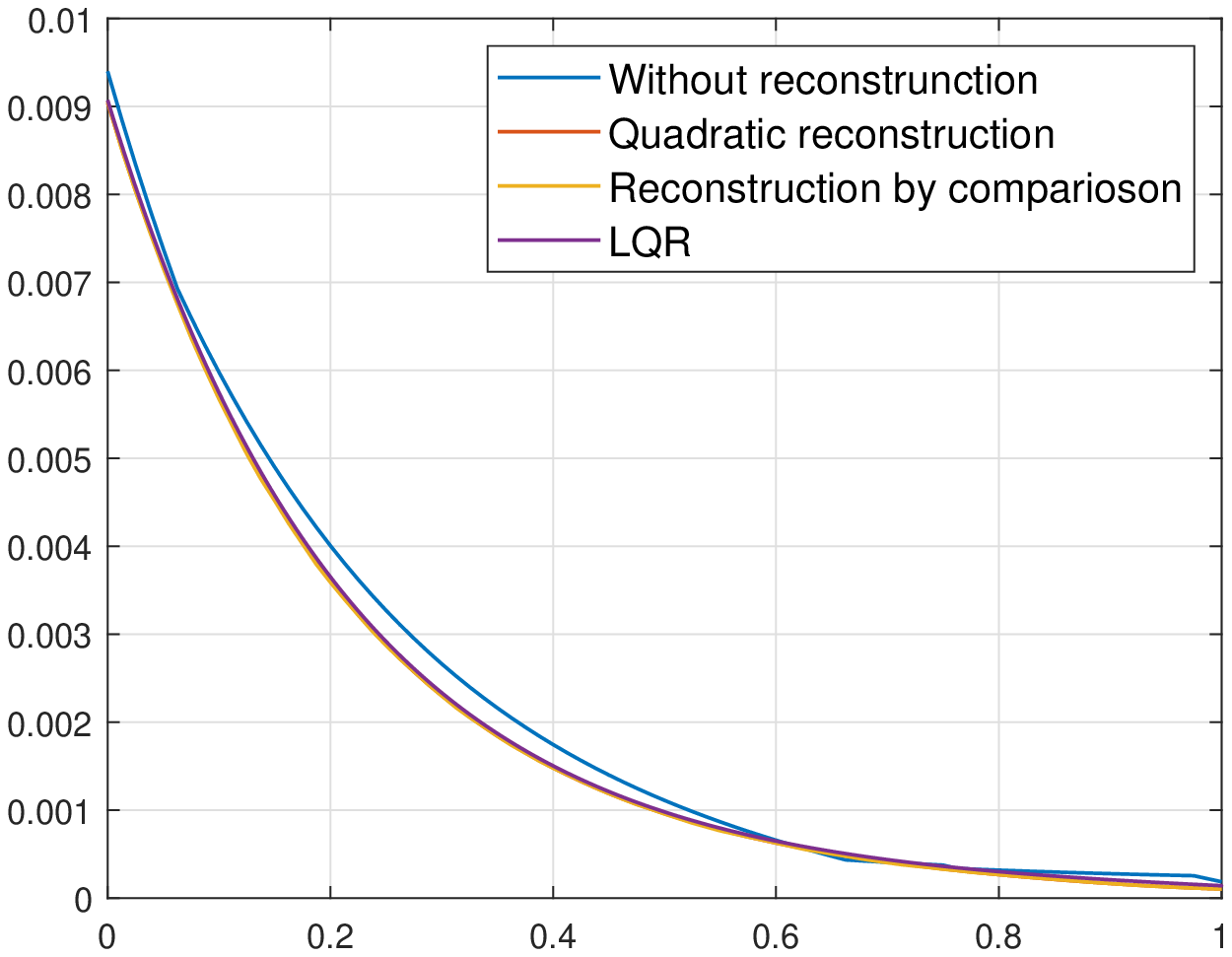}
       \includegraphics[scale=0.45]{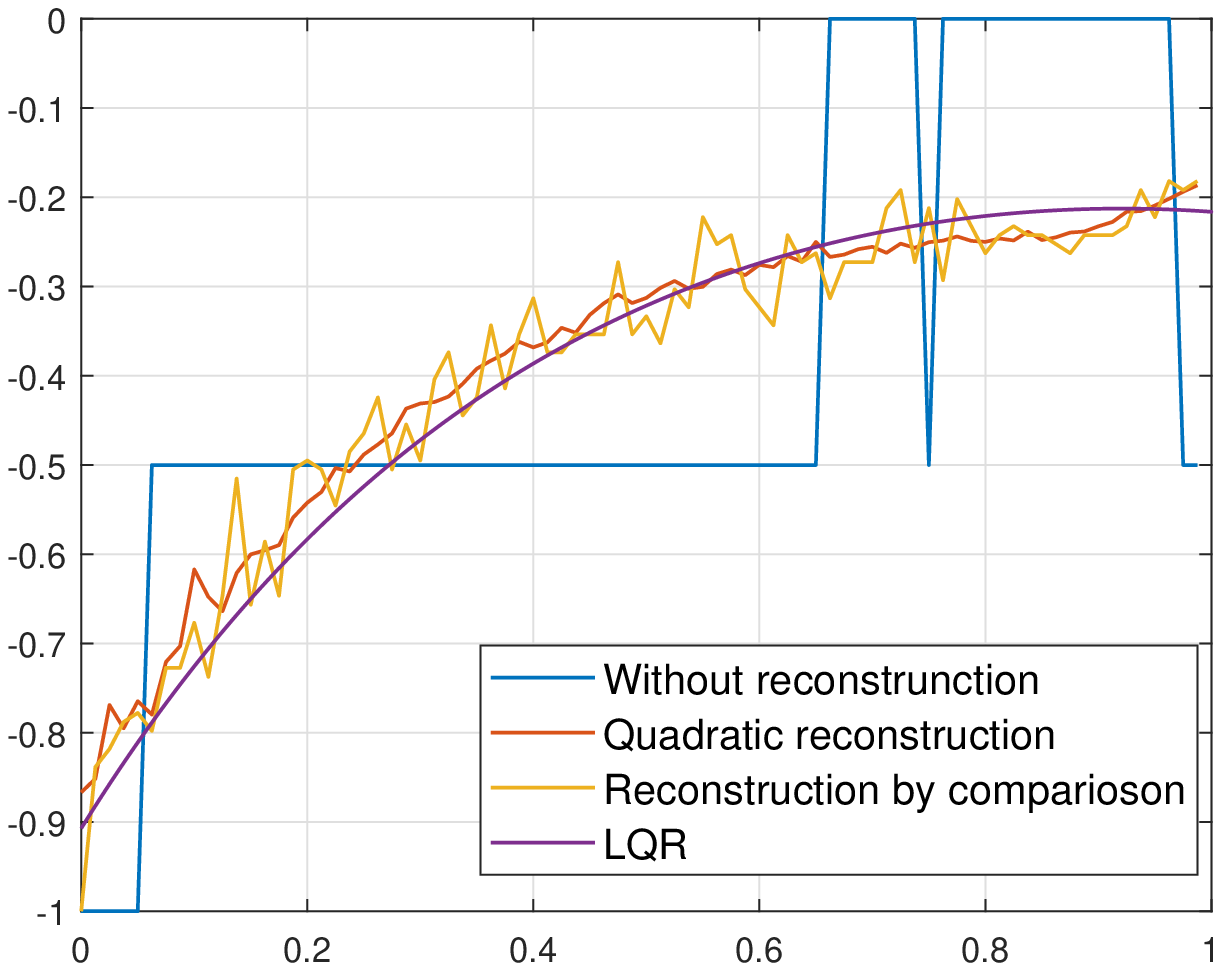} 
       \caption{Test 1: Cost functional (left) and optimal control (right) with different techniques for the feedback reconstruction. }
       \label{fig3:cost11feed}
	\end{figure}

In Figure \ref{fig3:cost11feed} we show the cost functional and the optimal control with all these techniques. In particular, it is possible to see from the plot of the optimal control that the quadratic feedback reconstruction is more stable, while the reconstruction by comparison presents a scattering behaviour.
	\begin{table}[htbp]
		\centering
		\begin{tabular}{cccccccc}
					\toprule
					 $\Delta t$  & CPU & $Err_{2}$ &  $Err_{\infty} $ & $Order_{2}$  & $Order_{\infty}$  \\
					\midrule
					0.1 &    \phantom{x} 0.02s&   0.251 &   0.229 &   & \\
					0.05 &     \phantom{x} 0.04s &   0.109 &  9.5e-2 &  1.21 &  1.27 \\
					0.025 &    \phantom{x} 0.63s & 3.3e-2 & 3.0e-2  &   1.71 &    1.65 \\
					0.0125 &     \phantom{x}24s & 1.1e-2 &   5.9e-3 & 1.58 &   2.36  \\
			\bottomrule		
		\end{tabular}
   	\caption{Test 1: Error analysis and order of convergence for TSA-POD method with $\ep= \Delta t^2$, $2$ POD basis and quadratic reconstruction.}
	\label{test1:2contrpar}
	\end{table}

\section{Conclusions}

In this work we have presented two algorithms to reconstruct the control in feedback form based on a tree structure proposed in \cite{AFS18}. Given the knowledge of the value function on the nodes of the tree, it is possible to introduce interpolation operators on scattered data to obtain the synthesis of the feedback including more discrete controls. This technique has been coupled with a Model Order Reduction method to reduce the dimensionality of the problem, allowing a fast and accurate computation of the optimal control. In the present work we have restricted ourselves to low dimensional reconstructions. In the next future our aim is to extend this idea to more high dimensional general problem using e.g. kernel methods \cite{F15}.

\bibliography{bib}

\end{document}